\documentclass[10pt,twoside]{article}
\usepackage{Latex-document}
\usepackage{amsmath,amssymb,amscd,enumerate}
\usepackage{Latex-document}
\usepackage{amsmath,amssymb}
\def\volumeno{I}
\title{\bf The Power Set Function\thanks{Partially supported by the
Isreal Science Foundation.}\vskip 6mm}
\author{Moti Gitik\thanks{School of Mathematical Sciences,
Tel Aviv University, Tel Aviv, Israel. E-mail: \ \ \ \ \ \ \ \ \ \
\ \ \ \ \ \ \ \ \ \ \ \ \ \ \ \ \ \
gitik@post.tau.ac.il}\vspace*{-0.5cm}}
\date{\vspace{-8mm}}

\begin{document}

\maketitle

\thispagestyle{first} \setcounter{page}{507}

\begin{abstract} \vskip 3mm

We survey  old and recent results on the problem of finding a complete set of rules describing
the behavior of the power function, i.e. the function which takes a cardinal
$\kappa$ to the cardinality of its power $2^\kappa$.
\end{abstract}

\vskip 12mm

\section{Introduction}
\vskip-5mm \hspace{5mm}

    One of the central topics of Set Theory since Cantor was the
study of the power function . The basic
problem is to determine all the possible values of $2^\kappa$ for a
cardinal $\kappa$. Paul Cohen \cite{Coh} proved the independence of the
Continuum Hypothesis and invented the method of forcing. Shortly
after, Easton \cite{Eas} building on Cohen's results showed  the function
$\kappa \longrightarrow 2^\kappa$, for regular $\kappa$, can behave in any prescribed
way consistent with K\"onig's Theorem. This reduces
the study to singular cardinals. It turned out that the situation with powers of singular
cardinals is much more involved. Thus, for example, a remarkable theorem of Silver \cite{Sil}
states that a singular cardinal of uncountable cofinality cannot be first to violate GCH.
The Singular Cardinal Problem is the problem of finding a complete
set of rules describing the behavior of the power function on singular cardinals.
 There are three main tools for dealing with the problem:
$pcf$-theory, inner models theory and forcing involving large cardinals.

\section{Classical results and basic definitions}
\vskip-5mm \hspace{5mm}

       In 1938 G\"odel proved the consistency of the Axiom of Choice (AC) and the Generalized
       Continuum Hypothesis (GCH) with the rest axioms of set theory. In 1963 Cohen proved the
        independence of AC and GCH.
He showed, in particular, that $2^{\aleph_0}$ can be arbitrary
large. Shortly after Solovay proved that $2^{\aleph_0}$ can take
any value $\lambda$ with
 $cf (\lambda) > \aleph_0$.
The $cofinality$ of a limit ordinal $\alpha$ ($cf (\alpha)$) is the
least ordinal $\beta \leq \alpha$ so that there is a function $f : \beta \longrightarrow \alpha$
 with $rng (f)$ unbounded in $\alpha$.
A cardinal $\kappa$ is called a $regular$ if  $\kappa = cf(\kappa)$.
Otherwise a cardinal is called a $singular$ cardinal. Thus, for example, $\aleph_8$ is regular
 and $\aleph_\omega$ is singular of cofinality $\omega$.

By a result of Easton \cite{Eas}, if we restrict ourselves to regular
cardinals, then every class function $F : Regulars \longrightarrow
Cardinals$ $satisfying$

(a) $\kappa \leq \lambda$  $implies$ $F(\kappa) \leq F(\lambda)$

(b) $cf ( F(\kappa) > \kappa$ ( K\"onig's Theorem)

\noindent can be realized as a power function in a generic extension.

From this point we restrict ourselves to singular cardinals.

\section{Restrictions on the power of singular cardinals}
\vskip-5mm \hspace{5mm}

The Singular Cardinal Problem (SCP) is the problem of finding a complete
set of rules describing the behavior of the power function on singular
cardinals. For singular cardinals there are more limitations. Thus

(c) (Bukovsky - Hechler) If $\kappa$ is a singular and there is $\gamma_
0 < \kappa$ such that $ 2^ \gamma = 2 ^{\gamma_ 0}$  for every $\gamma,
\gamma\leq\gamma<\kappa$ , then $2^ \kappa = 2^{\gamma_ 0}$.

(d) (Silver)  If  $\kappa$ is a singular strong limit cardinal of
uncountable cofinality and  $2^\kappa > \kappa ^+$  then $\{ \alpha <
\kappa  | 2^\alpha  > \alpha^ + \}$ contains a closed unbounded subset
of $\kappa$.

A set $C \subset \kappa$ is called a $closed$ $unbounded$ subset of
$\kappa$ iff

(1) $\forall \alpha < \kappa \exists \beta \in C (\beta > \alpha)$ ($unbounded$)

(2) $\forall \alpha < \kappa ( C \cap \alpha \not = \phi \Rightarrow sup ( C \cap \alpha) \in C)$  ($closed$).

Subsets of $\kappa$ containing a closed unbounded set form a filter over
$\kappa$ which is $\kappa$ complete. A positive for this filter sets are
called $stationry$.

(e) (Galvin - Hajnal, Shelah) If  $\aleph_\delta$ is strong limit and
$\delta< \aleph_\delta$   then $2^{\aleph_ \delta} <
\aleph_{2^{|\delta|^+}}$

(f) (Shelah) It is possible to replace $2^{|\delta|^+}$ in (e) by $ |\delta|^{+4}$.

(g) (Shelah) Let $\aleph_\delta$ be the $\omega_1$ -th fixed point
of the $\aleph$ - function. If it is a strong limit , then
$2^{\aleph_\delta} < min( ( 2^{\omega_1})^+$ -fixed point,
$\omega_4$ -th fixed point).

A cardinal $\kappa$ is called a $fixed$ $point$ of the $\aleph$ function
if $\kappa = \aleph_\kappa$.

It is unknown if $4$ in (f) and in (g) can be reduced or just replaced
by $1$. One of the major questions in Cardinal Arithmetic asks if
$2^{\aleph_\omega}$ can be bigger than $\aleph_{\omega_1}$ provided it
is a strong limit. We refer to the books by Jech \cite{Jec} and by
Shelah \cite{She} for the proofs of the above results.

\section{Inner models and large cardinals}
\vskip-5mm \hspace{5mm}

There are other restrictions which depend on large cardinals. Thus
the celebrated Covering Theorem of  Jensen \cite{Dev} implies that
for every singular strong limit cardinal $\kappa$   $2^\kappa =
\kappa^+$, provided  the universe is close to G\"odel's model $L$
( precisely, if $o \#$ does not exist, or, equivalently , there is
no elementary embedding from $L$ into $ L$). On the other hand,
using large cardinals (initially supercompact cardinals were used
\cite{Kan}) it is possible to have the following.

(Prikry-Silver, see \cite{Jec}):

$\kappa$ is a strong limit of cofinality $\omega$ and $2^\kappa >
\kappa^+$.

(Magidor \cite{Mag1},\cite{Mag2},\cite{Mag3}):

(1) the same with $\kappa$ of any uncountable cofinality.

(2) the same with $\kappa = \aleph _ \omega$.

So, the answer to SCP may depend on presence of particular large
cardinals. Hence, it is reasonable to study the possibilities for
the power function level by level according to existence of
particular large cardinals. There are generalizations of the
G\"odel model $L$ which may include bigger and bigger large
cardinals, have nice combinatorial properties,  satisfy GCH and
are invariant under set forcing extensions. This models are called
$Core$ $Models$. We refer to the book by Zeman \cite{Zem} for a
recent account on this fundamental results.

The Singular Cardinals Problem can now be reformulated as follows:

Given a core model $K$ with certain large cardinals. Which functions can
be realized in extensions of $K$ as power set functions , i.e. let $F :
Ord \longrightarrow Ord$ be a class function in $K$, is there an
extension (generic) of $K$ satisfying $2^{\aleph_\alpha} = \aleph
_{F(\alpha)}$ for all ordinals $\alpha$?

We will need few definitions.

An uncountable cardinal $\kappa$ is called a $measurable$ cardinal iff
there is $\mu : P(\kappa) \longrightarrow$ $\{0,1\}$ such that

(1) $\forall \alpha < \kappa$   $\mu(\{\alpha\}) = 0$.

(2) $\mu(\kappa) = 1$.

(3) $A \subseteq B \Longrightarrow \mu(A) \leq \mu(B)$.

(4) $\forall \delta < \kappa$ $ \forall$ $\{A_\nu | \nu < \delta \}$
subsets of $\kappa$ with $\mu(A_\nu) = 0$  $ \mu (\cup$ $A_\nu | \nu
<\delta\})= 0$.

If $\kappa$ is a measurable, then it is possible always to find $\mu$
with an additional property  called $normality$:

(5) If $\mu(A) = 1$ and $f : A \longrightarrow \kappa,  f(\alpha) <
\alpha$ then there is a subset of $A$ of measure one on which $f$ is
constant. Further by measure we shall mean a normal measure, i.e. one
satisfying (1)--(5). A cardinal $\kappa$  has the $Mitchell$ $order$
$\geq 1 (o(\kappa)\geq1)$ iff $ \kappa$ is a measurable. A cardinal
$\kappa$ has the $Mitchell$ $order$ $\geq 2$ ($o(\kappa)\geq 2$) iff
there is a measure over $\kappa$ concentrating on measurable cardinals,
 i.e.
$\mu(\{\alpha< \kappa | o(\alpha)\geq1 \} ) = 1.$

In a similar fashion we can continue further , but up to $\kappa^{++}$
only. Just the total number of ultrafilters over $\kappa$ under GCH is
$\kappa^{++}$. In order to continue above this point , directed systems
of ultrafilters called extenders are used. This way we can reach
$\kappa$ with $o(\kappa) = Ord$. Such $\kappa$ is called a $strong$
$cardinal$. Core models are well developed to the level of strong
cardinal and much further. Almost all known consistency results on the
Singular Cardinals Problem require large cardinals below the level of  a
strong cardinal.

\section{Finite gaps}
\vskip-5mm \hspace{5mm}

By results of Jensen \cite{Dev}, Dodd- Jensen  \cite{J-D}, Mitchell
\cite{Mit}, Shelah \cite{She} and Gitik \cite{Gitsch} nothing
interesting in sense of SCP happens bellow the level of $o(\kappa) =
\kappa ^ {++}$. If there is $n < \omega$ such that for every $\alpha,
o(\alpha) \leq \alpha ^{+n}$, then we have the following additional
restrictions:

(1) (Gitik-Mitchell \cite{Git-Mit} ) If  $\kappa$ is a singular
strong limit and $2^\kappa = \kappa ^{+m}$ for some $m > 1$, then,
in $K, o(\kappa) \geq \kappa^{+m}$. In particular, $m \leq n$.

(2) If $\kappa$ is a singular cardinal of uncountable cofinality and for
some $m, 1\leq m <\omega$  $\{\alpha < \kappa | 2^\alpha = \alpha
^{+m}\}$ is stationry, then $\{\alpha < \kappa | 2^\alpha = \alpha
^{+m}\}$ contains a closed unbounded subset of $\kappa$.

By results of Merimovich \cite{Mer} it looks like this are the only
restrictions.

\section{Uncountable cofinality case}
\vskip-5mm \hspace{5mm}

Assume only that there is no inner model with a strong cardinal. Then we
have the following restrictions:

(1) If $\kappa$ is a singular strong limit cardinal of uncountable
cofinality $\delta$ and $2^\kappa \geq \lambda > \kappa ^+$ , where
$\lambda$ is not the successor of a cardinal of cofinality less than
$\kappa$,  then $o(\kappa) \geq \lambda + \delta$, if $\delta >
\omega_1$ or $o(\kappa) \geq \lambda$, if $\delta = \omega _1$.

(2) Let $\kappa$ be a singular strong limit cardinal of uncountable
cofinality $\delta$ and let $\tau < \delta$. If $A= \{\alpha < \kappa |
cf \alpha > \omega , 2^\alpha = \alpha ^{+\tau}\}$ is stationry, then
$A$ contains a closed unbounded subset of $\kappa$.

(3) If $\delta < \aleph_\delta , \aleph_\delta$ strong limit then
$2^{\aleph_\delta} < \aleph_{ | \delta|^+ }$

(This was improved recently by R.Schindler [11] to many Woodin
cardinals).

(4) Let $\aleph_\delta$ be the $\omega_1$ -th fixed point of the $\aleph
$-function. If it is a strong limit cardinal then $2^{\aleph_\delta} <
\omega _2$ -th fixed point.

(5) If $a$ is an uncountable set of regular cardinals with $min( a)
>2^{|a|^++\aleph_2}$, then $| pcf(a)| = | a|$, where $pcf(a) = \{cf (\Pi a/D) | D$ is an ultrafilter over $a\}.$

It is a major problem of Cardinal Arithmetic if it is possible to have a
set of regular cardinals $a$ with $min( a) >|a|$ such that $| pcf(a)| >
| a|$ . The results above were proved in Gitik-Mitchell \cite{Git-Mit},
and in \cite{Gitgch}. It is unknown if there is no further restrictions
in this case (i.e. singulars of uncountable cofinality under  the
assumption that there is no inner model with a strong cardinal). Some
local cases were  checked by Segal \cite{Seg} and Merimovich
\cite{Mer1}.

\section{Countable cofinality case}
\vskip-5mm \hspace{5mm}

In this section we revue some more recent results dealing with countable cofinality.
First suppose that

$(\forall n<\omega \exists \alpha$   $o(\alpha)= \alpha^{+n})$,
but $\neg (\exists \alpha$  $o(\alpha) = \alpha ^{+\omega})$.

Then the following holds: Let $\kappa$ be a cardinal of countable
cofinality such that for every $n< \omega$  $\{\alpha < \kappa |
o(\alpha) \geq \alpha ^{+n}\}$ is unbounded in $\kappa$. Then for
every $\lambda \geq \kappa ^ +$ there is a cardinal preserving
generic extension satisfying ``$\kappa$ is a strong limit and
$2^\kappa \geq \lambda$ ". So the gap between a singular and its
power can already be arbitrary large .  But by \cite{Gitgch}:

If $2^\kappa \geq \kappa ^{+\delta}$ for $\delta \geq \omega _1 $, then
GCH cannot hold below $\kappa$.

(Actually , GCH can hold if the gap is at most countable
\cite{Gitbl}.)

We do not know if ``pcf (a)uncountable for a countable $a$" is
stronger than the assumption above. If  we require also GCH below,
then it is.

Once one likes to have  uncountable gaps between a singular
cardinal and its power together with GCH below , then the
following results provide this and are sharp. The proofs are
spread through papers \cite{Gitbl}, \cite{Git-Mag},
\cite{Git-Mit}.

Suppose that $\kappa>\delta\geq\aleph_0, \delta$ is a cardinal, $2^\kappa \geq \kappa^{+\delta}, cf \kappa = \aleph_0$ and GCH below $\kappa$. Then

(i) $cf \delta = \aleph_0$ implies (that in the core model) for every $\tau < \delta$  $\{\alpha< \kappa | o(\alpha) \geq \alpha^{+\tau}\}$ is unbounded in $\kappa$.

(ii) $cf \delta > \aleph_0$ implies (in the core model) $o(\kappa) \geq \kappa^{+\delta + 1} +1$  or  $\{\alpha< \kappa | o(\alpha)\geq \alpha^{+\delta + 1 }\}$ is unbounded in $\kappa$.

Finally let us consider the following large cardinal: $\kappa$ is
singular of cofinality $\omega$ and for every $\tau < \kappa$ $\{\alpha<
\kappa | o(\alpha)\geq \alpha^{+\tau}\}$ is unbounded in $\kappa$.

Under this assumption it is possible to blow up the power of
$\kappa$ arbitrary high preserving GCH below $\kappa$. Also, it is
possible to turn $\kappa$ into the first fixed point of the
$\aleph$ function, see \cite{Gitffp}. This answers Question
($\gamma$) from the Shelah's book on cardinal arithmetic
\cite{She}. What are the possibilities for the power function
under the assumption above? First in order to be able to deal with
cardinals above $\kappa$, let us replace it by a global one:

For every $\tau$ there is $\alpha$  $o(\alpha)\geq \alpha^{+\tau}$.

We do not know the status of ``$pcf$ of a countable set uncountable'',
but other limitations like

(1) $\aleph_\omega$ strong limit implies $2^{\aleph_\omega} < \aleph_{\omega_1}$

(2) If $\kappa$ is a singular of uncountable cofinality then either
$\{\alpha< \kappa | 2^\alpha\geq \alpha^+\}$ $ or$  $\{\alpha< \kappa |
2^\alpha> \alpha^+\}$ contains a closed unbounded subset of $\kappa$

\noindent are true below strong cardinal.

By recent result \cite{G-Sh-Sc} the negation of the second assumption
implies initially unrelated statement -  Projective Determinacy. We
refer to the books by A. Kanamori \cite{Kan} and H. Woodin \cite{Woo} on
this subject. We conjecture that there is no other limitations, i.e. (1)
with $\aleph _\omega$ replaced by $\aleph _\delta$ for $\delta<
\aleph_\delta$, (2) and the classical ones.

\section{One idea}
\vskip-5mm \hspace{5mm}

Let us conclude with a sketch of  one basic idea which is crucial for
the forcing constructions in the countable cofinality case. Let $U$ be a
$\kappa$ complete nontrivial ultrafilter over $ \kappa$ (say, in $K$). A
sequence $\langle \delta_n | n<\omega \rangle$ is called a $Prikry$
$sequence$ for $U$ iff for each $A \in U$ $ \exists n_0 \forall n \geq
n_0$  $ \delta_n \in A$. Suppose now that $\kappa$ is a strong limit
singular cardinal of cofinality $\omega$ and $2^\kappa = \kappa^{++}$.
Then, usually (by \cite{Gitsch}, \cite{Git-Mit}), we will have a
sequence $\langle U_\alpha | \alpha < \kappa^{++}\rangle$ of
ultrafilters in $K$ and a sequence $\langle \delta_{\alpha,n} | \alpha <
\kappa^{++} , n < \omega\rangle$ so that

(1) $\alpha< \beta \Longrightarrow \exists n_0 \forall n \geq n_0
\delta_{\alpha,n} < \delta_{\beta,n}$,

(2) $\langle \delta_{\alpha,n} | n< \omega\rangle$ is a Prikry sequence for $U_\alpha$.

Ultrafilters $U_\alpha$ are different here. So each sequence $\langle
\delta_{\alpha,n} | n< \omega\rangle $  relates to unique ultrafilter
from the list. But once $\kappa^{++}$ is replaced by $\kappa^{+++}$, the
corresponding sequence of ultrafilters  $\langle U_\alpha |\alpha <
\kappa^{+++}\rangle$ will have different $\alpha$ and $\beta$,
$\kappa^{++} <\alpha< \beta < \kappa^{+++}$  with $U_\alpha = U_\beta$.
Then a certain Prikry sequence $\langle \delta_n | n< \omega \rangle$
may pretend to correspond to both $U_\alpha$ and $U_\beta$. In order to
decide, we will need a Prikry sequence for some $U_\gamma$ with $\gamma
< \kappa^{++}$ (more precisely,  if $f_\beta$ is the canonical one to
one correspondence in $K$ between $\kappa^{++}$ and $\beta$ then
$f_\beta(\gamma) = \alpha$). Dealing with $\kappa^{+4}$ we will need go
down twice, first to $\kappa^{+3}$ and after that to $\kappa^{++}$. In
general, for $n, 3\leq n<\omega$, $n-2$-many times. Certainly, it is
impossible to go down infinitely many times, but instead we replace the
fixed $\kappa$ by an increasing sequence $\langle \kappa_n | n<
\omega\rangle$ with each $\kappa_n$ carrying  $\kappa_n^{+n+3}$ many
ultrafilters. Now it turns out to be possible to add $\omega$ -
sequences with no assignment to ultrafilters. Just the number of steps
needed to produce the assignment is $\omega$ which is not enough for
sequences of the length $\omega$.

\label{lastpage}

\end{document}